\documentclass[amsbook,12pt]{article}
\usepackage{epsfig, amsfonts,amssymb, amsmath}
\usepackage[usenames]{color}
\usepackage{picinpar}
\newcommand{\nl}{\mbox{}\\}

\begin{document}
%
% \thispagestyle{empty}
%
% ----------------------------------------------------------------
%
\mbox{} \vspace{-2.000cm} \\
\begin{center}
{\Large \bf %
Some remarks on the regularity time of Leray} \\
\mbox{} \vspace{-0.250cm} \\
\mbox{\Large \bf %
solutions to the Navier-Stokes equations} \\
\nl
\mbox{} \vspace{-0.200cm} \\
{\large \sc P\!\:\!.\;Braz e Silva,}$\mbox{}^{\!\!\;\!1}$
{\large \sc J.\;P\!\:\!.\;Zingano}$\mbox{}^{\:\!2}$
{\large \sc and P\!\:\!.\;R.\;Zingano}$\mbox{}^{\:\!2}$ \\
\mbox{} \vspace{-0.150cm} \\
{\small $\mbox{}^{1}\:\!$Departmento de Matem\'atica} \\
\mbox{} \vspace{-0.670cm} \\
{\small Universidade Federal de Pernambuco} \\
\mbox{} \vspace{-0.670cm} \\
{\small Recife, PE 50740, Brazil} \\
\mbox{} \vspace{-0.300cm} \\
$\mbox{}^{2}\:\!${\small Departamento de Matem\'atica Pura e Aplicada} \\
\mbox{} \vspace{-0.670cm} \\
{\small Universidade Federal do Rio Grande do Sul} \\
\mbox{} \vspace{-0.670cm} \\
{\small Porto Alegre, RS 91509, Brazil} \\
\nl
\mbox{} \vspace{-0.300cm} \\
{\bf Abstract} \\
\mbox{} \vspace{-0.525cm} \\
%
% \begin{minipage}[t]{12.125cm}
% \begin{minipage}[t]{11.250cm}
%
\begin{minipage}[t]{11.000cm}
{\small
\mbox{} \hspace{+0.200cm}
We strengthen the classic result
about the regularity time $t_{\ast} $
of arbitrary Leray solutions
of the Navier-Stokes equations
in $ \mathbb{R}^{n} $ ($n = 3, 4$),
which have the form
$ \;\!t_{\ast} \!\:\!\leq\!\;\! K_{\mbox{}_{\!\;\!3}} \,\nu^{-\,5}\:\!
\|\, \mbox{\boldmath $u$}_0 \;\!
\|_{L^{2}(\mathbb{R}^{3})}^{\:\!4} $
if $ \:\!n = 3 $,
and
$ \!\;\!\:\!t_{\ast} \!\:\!\leq\!\;\! K_{\mbox{}_{\!\;\!4}}
\,\nu^{-\,3}\:\!
\|\, \mbox{\boldmath $u$}_0 \;\!
\|_{L^{2}(\mathbb{R}^{4})}^{\:\!2} \!\;\!$
if $ \:\!n = 4 $
(in particular, by
reducing \linebreak
the current known values
for the constants
$ \:\!K_{\mbox{}_{\!\:\!3}} \!\;\! $,
$ \!\;\!K_{\mbox{}_{\!\:\!4}} $).
Some related results
are also included
in our discussion. \\
}
\end{minipage}
\end{center}
%
%
% -------------------------------------------------------
%
\nl
\setcounter{page}{1}
\mbox{} \vspace{-0.800cm} \\
%
%
% -----------------------------------------------------
%
%
% *************************************************************
% *                                                           *
% *                 Section 1: Introduction                   *
% *                                                           *
% *************************************************************
%

{\bf 1. Introduction} \\
\mbox{} \vspace{-0.750cm} \\

In this note
we rederive and slightly improve
a fundamental result
originally obtained by J.\;Leray \cite{Leray1934}
in dimension $ n = 3 $
(see (1.3) below)
for the global weak solutions
introduced in \cite{Leray1934}
to solve the incompressible
Navier-Stokes equations \\
\mbox{} \vspace{-0.650cm} \\
\begin{equation}
\tag{1.1$a$}
\mbox{\boldmath $u$}_{t}
\;\!+\,
\mbox{\boldmath $u$} \!\cdot\! \nabla \mbox{\boldmath $u$}
\,+\,
\nabla p
\;=\;
\nu \,
\Delta\:\! \mbox{\boldmath $u$},
%
% \qquad
% \nabla \!\cdot  \mbox{\boldmath $u$}(\cdot,t) \,=\, 0,
%
\end{equation}
\mbox{} \vspace{-0.975cm} \\
\begin{equation}
\tag{1.1$b$}
\nabla \!\cdot  \mbox{\boldmath $u$}(\cdot,t) \,=\, 0,
\end{equation}
\mbox{} \vspace{-0.975cm} \\
\begin{equation}
\tag{1.1$c$}
\mbox{\boldmath $u$}(\cdot,0)
\,=\,
\mbox{\boldmath $u$}_0 \in
\mbox{\boldmath $L$}^{2}_{\sigma}(\mathbb{R}^{n}),
\end{equation}
\mbox{} \vspace{-0.250cm} \\
and which is valid more generally
for arbitrary solutions
$ {\displaystyle
\;\!
\mbox{\boldmath $ u $}(\cdot,t)
\in
L^{\infty}((0,\infty), \:\!\mbox{\boldmath $L$}^{2}_{\sigma}(\mathbb{R}^{n}))
} $
$ {\displaystyle
\cap \,
L^{2}((0,\infty), \:\!\dot{\mbox{\boldmath $H$}}
\mbox{}^{\!\;\!1}\!\;\!(\mathbb{R}^{n}))
\cap \;\!
C_{\mbox{\scriptsize w}}([\;\!0, \infty),
\mbox{\boldmath $L$}^{2}(\mathbb{R}^{n}))
} $
of (1.1),
where $ n = 3 $ or $4$,
satisfying
the so-called
{\em strong energy inequality} \\
\mbox{} \vspace{-0.725cm} \\
\begin{equation}
\tag{1.2}
\|\, \mbox{\boldmath $u$}(\cdot,t) \,
\|_{\mbox{}_{\scriptstyle L^{2}(\mathbb{R}^{n})}}
  ^{\:\!2}
\!\:\!+\,
2 \, \nu \!\!\;\!
\int_{\mbox{}_{\mbox{\footnotesize \!\,\!$s$}}}
    ^{\mbox{\footnotesize \;\!$t$}}
\!\!\:\!
\|\, D \,\!\mbox{\boldmath $u$}(\cdot,\tau) \,
\|_{\mbox{}_{\scriptstyle L^{2}(\mathbb{R}^{n})}}
  ^{\:\!2}
d\tau
\,\leq\:
\|\, \mbox{\boldmath $u$}(\cdot,s) \,
\|_{\mbox{}_{\scriptstyle L^{2}(\mathbb{R}^{n})}}
  ^{\:\!2}
\!\:\!,
\quad
\forall \;\;\! t \geq s
\end{equation}
\mbox{} \vspace{-0.750cm} \\
for a.e.$\;s \geq 0 $,
including $ s = 0 $.\footnote{%
%
% ----------------------------------------------
%
%              Footnote 1
%
% ----------------------------------------------
%
For the definition of
$ {\displaystyle
\:\!
\|\, \mbox{\boldmath $u$}(\cdot,t) \,
\|_{\scriptstyle L^{2}(\mathbb{R}^{n})}
\!\;\!
} $,
$ {\displaystyle
\|\, D \,\!\mbox{\boldmath $u$}(\cdot,t) \,
\|_{\scriptstyle L^{2}(\mathbb{R}^{n})}
\!\:\!
} $
and other similar norms,
% of $ \mbox{\boldmath $u$}(\cdot,t) $,
see (1.5) next.
}
%
% ------------------------------- End of Footnote 1
%
Such solutions are now
called {\em Leray}
(or {\em Leray-Hopf\/})
solutions,
after \cite{Hopf1950, Leray1933, Leray1934}.
They have also been constructed
in higher dimensions
and other methods,
see e.g.\;\cite{% ConstantinFoias1988,
CaffarelliKohnNirenberg1982,
DoeringGibbon1995, FujitaKato1964,
Giga1986, Hopf1950, Kato1984, Ladyzhenskaya1969,
Lions1969, RobinsonRodrigoSadowski2016, Serrin1963,
ShinbrotKaniel1966, Sohr2001, Temam1984, Wahl1985}.
% and references therein.
% CaffarelliKohnNirenberg1982, %
% FujitaKato1964, % Galdi2000, Giga1986,
% GigaMiyakawa19xx, Heywood1980,
% Hopf1950, % KajikiyaMiyakawa1986,
% Kato1984, KiselevLadyzhenskaya1957, %
% Ladyzhenskaya1969,
% Lions1969, Masuda1984,
% Serrin1963,
% ShinbrotKaniel1966, Sohr2001, Temam1984,
% SohrWahl1987, Weissler1980}.
%
In (1.1),
$ \nu > 0 $ is a given constant,
\mbox{$ \mbox{\boldmath $u$} = \mbox{\boldmath $u$}(x,t) $}
and
$ p = p(x,t) $
are the unknowns
(the flow velocity and pressure, respectively),
and condition (1.1$c$)
is meant in $ \mbox{\boldmath $L$}^{2}(\mathbb{R}^{n}) $,
i.e.,
$ {\displaystyle
\,\!\,\!
\|\, \mbox{\boldmath $u$}(\cdot,t) - \mbox{\boldmath $u$}_0 \;\!
\|_{L^{2}(\mathbb{R}^{n})}
\!\;\!\rightarrow 0
\,\!\,\!\,\!
} $
as $ \,\!\,\! t \,\mbox{\footnotesize $\searrow$}\,0 $.
$\!\;\!$As usual,
$ \!\;\!\mbox{\boldmath $L$}^{2}_{\sigma}(\mathbb{R}^{n}) $
is the
space
of solenoidal fields
$ \:\!\mbox{\bf v} = (v_{1}, ..., v_{n}) \!\:\!\in
\mbox{\boldmath $L$}^{2}(\mathbb{R}^{n}) \equiv L^{2}(\mathbb{R}^{n})^{n} \!\:\!$
with
$ \nabla \!\cdot \mbox{\bf v} \!\;\!= 0 $
% \mbox{in $ {\cal D}^{\prime}(\mathbb{R}^{n}) $}
% (i.e., as distributions)
in the distributional sense, \linebreak
$ \dot{\mbox{\boldmath $H$}}\mbox{}^{1}(\mathbb{R}^{n}) =
\dot{H}^{1}(\mathbb{R}^{n})^{n} \!$,
where
$ \dot{H}^{1}(\mathbb{R}^{n}) $
denotes the homogeneous Sobolev space of order~1
(see e.g.\;\cite{BahouriCheminDanchin2011}, p.\;25),
and
$ \mbox{\boldmath $H$}\mbox{}^{m}(\mathbb{R}^{n}) =
H^{m}(\mathbb{R}^{n})^{n}\!$,
where
$ H^{m}(\mathbb{R}^{n}) $
is the Sobolev space of order $m$
(\cite{BahouriCheminDanchin2011}, p.\;38).
Also,
$ C_{\mbox{\scriptsize w}}(I, \:\!\mbox{\boldmath $L$}^{2}(\mathbb{R}^{n})) $
denotes the set of mappings from a
given interval~\mbox{$ I \subseteq \mathbb{R} $}
to
$ \mbox{\boldmath $L$}^{2}(\mathbb{R}^{n}) $
that are $L^{2}$-\;\!weakly continuous
at each $\:\! t \in I$. \\
\mbox{} \vspace{-0.900cm} \\

The result about Leray solutions
of (1.1)
that concerns us here
is the following.
In dimension $ n = 3 $ or $ 4 $,
it is known that
$ {\displaystyle
\mbox{\boldmath $u$}
\in
C^{\infty}(\mathbb{R}^{n} \!\times \!\:\![\, t_{\ast}, \infty))
} $
and \\
\mbox{} \vspace{-0.575cm} \\
\begin{equation}
\tag{1.3$a$}
\mbox{\boldmath $u$}(\cdot,t)
\in
C^{0}([\,t_{\ast}, \infty), \:\!
\mbox{\boldmath $H$}\mbox{}^{\!\:\!m}\!\;\!(\mathbb{R}^{n})),
\quad \;\,
\forall \;\;\! m \geq 0
\end{equation}
\mbox{} \vspace{-0.225cm} \\
for some
regularity time
$ \:\!t_{\ast} \!\geq 0 \:\!$
that satisfies \\
\mbox{} \vspace{-0.625cm} \\
\begin{equation}
\tag{1.3$b$}
t_{\ast}
% \,\leq\;
% \mbox{\small $ {\displaystyle \frac{1}{\;\!128\;\!\pi^{2}} }$}
\;\!\leq\, K_{\mbox{}_{\!\:\!3}}
\,
\nu^{-\,5} \;\!
\|\, \mbox{\boldmath $u$}_0 \;\!
\|_{\mbox{}_{\scriptstyle L^{2}(\mathbb{R}^{3})}}
  ^{\:\!4}
\quad \;\;
(\:\!\mbox{if }\;\! n = 3\:\!),
\end{equation}
\mbox{} \vspace{-0.750cm} \\
\begin{equation}
\tag{1.3$c$}
t_{\ast}
% \,\leq\;
% \mbox{\small $ {\displaystyle \frac{1}{\;\!128\;\!\pi^{2}} }$}
\;\!\leq\, K_{\mbox{}_{\!\:\!4}}
\,
\nu^{-\,3} \;\!
\|\, \mbox{\boldmath $u$}_0 \;\!
\|_{\mbox{}_{\scriptstyle L^{2}(\mathbb{R}^{4})}}
  ^{\:\!2}
\quad \;\;
(\:\!\mbox{if }\;\! n = 4\:\!),
\end{equation}
\mbox{} \vspace{-0.250cm} \\
with constants
$ K_{\mbox{}_{\!\:\!3}} \!\;\!$,
$ K_{\mbox{}_{\!\:\!4}} $
that are independent of $ \;\!\nu $,
$ \mbox{\boldmath $u$}_0 $
or the solution $ \mbox{\boldmath $u$} $
considered.
For example,
it follows from
(\cite{Heywood1994}, p.\:14;
\cite{Leray1934}, p.\:246;
\cite{LorenzZingano2003}, p.\:18)
that \\
\mbox{} \vspace{-0.600cm} \\
\begin{equation}
\notag
K_{\mbox{}_{\!\:\!3}}
\:\!\leq\:
\mbox{\small $ {\displaystyle \frac{1}{\,128\;\!\pi^{2}} }$}
\;\!<\,
\mbox{\small $0.000\,791\;\!572$},
\end{equation}
\mbox{} \vspace{-0.225cm} \\
which appears to be the best
(smallest)
estimate presently known
for the constant
$ \!\;\!K_{\mbox{}_{\!\:\!3}} \!\;\! $ \linebreak
in the literature.
In Section 2,
it is shown that \\
\mbox{} \vspace{-0.600cm} \\
\begin{equation}
\tag{1.4}
K_{\mbox{}_{\!\:\!3}}
\:\!<\:
\mbox{\small $0.000\,464\,504\,284$},
\quad \;\;\;\,
K_{\mbox{}_{\!\:\!4}}
<\:
\mbox{\small $0.002\,727\,993\,110$},
\end{equation}
\mbox{} \vspace{-0.225cm} \\
along with other improvements
regarding the monotonic behavior
of $ \;\!\|\;\! D \mbox{\boldmath $u$}(\cdot,t) \,
\|_{\mbox{}_{\scriptstyle L^{2}(\mathbb{R}^{n})}} $
for large $ \!\;\!\:\!t $.
\!Moreover,
by a similar argument,
we note in {\small \sc Remark 2.3}
that,
for each $m$,
$ \|\, D^{m} \mbox{\boldmath $u$}(\cdot,t) \,
\|_{\mbox{}_{\scriptstyle L^{2}(\mathbb{R}^{n})}} \!\:\!$
also becomes monotonically decreasing for $\;\!t \gg 1$.
%
% with
% $ \;\!t^{\:\!m/2} \;\!\|\, D^{m} \mbox{\boldmath $u$}(\cdot,t) \,
% \|_{\mbox{}_{\scriptstyle L^{2}(\mathbb{R}^{n})}}
% \!\rightarrow \:\!0 $
% as $ \;\!t \rightarrow \infty $.
%
\!Our analysis is
inspired by the interesting approach to these questions
developed in \cite{KreissHagstromLorenzZingano2002, %
KreissHagstromLorenzZingano2003, Zhou2007, Zingano2015}.

%
% ---------------------------------------------------
%
%                     Notation
%
% ---------------------------------------------------
%
\mbox{} \vspace{-0.750cm} \\
{\bf Notation.}
As shown above,
boldface letters
are used for
vector quantities,
as in
$ {\displaystyle
\;\!
\mbox{\boldmath $u$}(x,t)
=
} $
$ {\displaystyle
(\:\! u_{\mbox{}_{\!\:\!1}}\!\;\!(x,t), ...\,\!\,\!,
 \:\! u_{\mbox{}_{\scriptstyle \!\;\! n}}\!\;\!(x,t) \:\!)
} $.
$\!$Also,
$ \nabla p \;\!\equiv \nabla p(\cdot,t) $
denotes the spatial gradient of $ \;\!p(\cdot,t) $,
$ D_{\!\;\!j} \!\;\!=\:\! \partial / \partial x_{\!\;\!j} $,
$ {\displaystyle
\,\!
\nabla \!\cdot \mbox{\boldmath $u$}
\:\!=
  D_{\mbox{}_{\!\:\!1}} u_{\mbox{}_{\!\:\!1}} \!\;\!+
  ... \!\;\!+
  D_{\mbox{}_{\scriptstyle \!\;\!n}} \,\!
  u_{\mbox{}_{\scriptstyle \!\;\!n}}
} $
is the (spatial) divergence of
$ \:\!\mbox{\boldmath $u$}(\cdot,t) $.
$ |\,\!\cdot\,\!|_{\mbox{}_{2}} \!\,\!\,\!$
denotes the Euclidean norm
in $ \mathbb{R}^{n} \!$,
and
$ {\displaystyle
\,\!\,\!
\| \:\!\cdot\:\!
\|_{\scriptstyle L^{q}(\mathbb{R}^{n})}
\!\;\!
} $,
$ 1 \leq q \leq \infty $,
are the standard norms
of the Lebesgue spaces
$ L^{q}(\mathbb{R}^{n}) $,
with the vector counterparts \\
\mbox{} \vspace{-0.625cm} \\
\begin{equation}
\tag{1.5$a$}
\|\, \mbox{\boldmath $u$}(\cdot,t) \,
\|_{\mbox{}_{\scriptstyle L^{q}(\mathbb{R}^{n})}}
\;\!=\;
\Bigl\{\,
\sum_{i\,=\,1}^{n} \int_{\mathbb{R}^{n}} \!
|\:u_{i}(x,t)\,|^{q} \,dx
\,\Bigr\}^{\!\!\:\!1/q}
\end{equation}
\mbox{} \vspace{-0.700cm} \\
\begin{equation}
\tag{1.5$b$}
\|\, D \mbox{\boldmath $u$}(\cdot,t) \,
\|_{\mbox{}_{\scriptstyle L^{q}(\mathbb{R}^{n})}}
\;\!=\;
\Bigl\{\,
\sum_{i, \,j \,=\,1}^{n} \int_{\mathbb{R}^{n}} \!
|\, D_{\!\;\!j} \;\!u_{i}(x,t)\,|^{q} \,dx
\,\Bigr\}^{\!\!\:\!1/q}
\end{equation}
\mbox{} \vspace{-0.350cm} \\
and, in general, \\
\mbox{} \vspace{-0.750cm} \\
\begin{equation}
\tag{1.5$c$}
\|\, D^{m} \mbox{\boldmath $u$}(\cdot,t) \,
\|_{\mbox{}_{\scriptstyle L^{q}(\mathbb{R}^{n})}}
\;\!=\;
\Bigl\{\!\!
\sum_{\mbox{} \;\;i, \,j_{\mbox{}_{1}} \!,..., \,j_{\mbox{}_{m}} =\,1}^{n}
\!\;\! \int_{\mathbb{R}^{n}} \!
|\, D_{\!\;\!j_{\mbox{}_{1}}}
\!\!\!\;\!\cdot\!\,\!\cdot\!\,\!\cdot \!\:\!
D_{\!\;\!j_{\mbox{}_{m}}}
\!\:\! u_{i}(x,t)\,|^{q} \,dx
\,\Bigr\}^{\!\!\:\!1/q}
\end{equation}
\mbox{} \vspace{-0.175cm} \\
if $ 1 \leq q < \infty $\/;
if $\, q = \infty $,
then
$ {\displaystyle
\;\!
\|\, \mbox{\boldmath $u$}(\cdot,t) \,
\|_{\mbox{}_{\scriptstyle L^{\infty}(\mathbb{R}^{n})}}
\!=\;\!
\max \, \bigl\{\,
\|\,u_{i}(\cdot,t)\,
\|_{\mbox{}_{\scriptstyle L^{\infty}(\mathbb{R}^{n})}}
\!\!: \, 1 \leq i \leq n
\,\bigr\}
} $, \linebreak
\mbox{} \vspace{-0.530cm} \\
$ {\displaystyle
\|\, D \,\!\mbox{\boldmath $u$}(\cdot,t) \,
\|_{\mbox{}_{\scriptstyle L^{\infty}(\mathbb{R}^{n})}}
\!=\;\!
\max \, \bigl\{\,
\|\, D_{\!\;\!j} \;\! u_{i}(\cdot,t)\,
\|_{\mbox{}_{\scriptstyle L^{\infty}(\mathbb{R}^{n})}}
\!\!: \:\! 1 \leq i, \:\!j \leq n
\,\bigr\}
} $
and,
for general \mbox{$m \!\;\!\geq\!\:\! 1$\/:} \\
\mbox{} \vspace{-0.500cm} \\
\begin{equation}
\tag{1.5$d$}
\|\, D^{m} \mbox{\boldmath $u$}(\cdot,t) \,
\|_{\mbox{}_{\scriptstyle L^{\infty}(\mathbb{R}^{n})}}
\!\;\!=\;
\max\,\Bigl\{\;\!
\|\;\! D_{\!\;\!j_{\mbox{}_{1}}}
\!\!\!\;\!\;\!\,\!\cdot \!\;\!\cdot \!\;\!\cdot \!\;\!\;\!
D_{\!\;\!j_{\mbox{}_{m}}}
\!\!\;\!\;\! u_{i}(\cdot,t)\,
\|_{\mbox{}_{\scriptstyle L^{\infty}(\mathbb{R}^{n})}}
\!\!\!\;\!: \;\!
1 \leq \!\;\!\;\!
i, \!\;\!\;\!j_{\mbox{}_{1}}\!\!\;\!\;\!, \!...\!\;\!\;\!,j_{\mbox{}_{m}}
\!\leq n
\!\;\!\;\!\Bigr\}.
\end{equation}
\mbox{} \vspace{-0.150cm} \\
The definitions chosen in (1.5) are very convenient
for the discussion that follows. \linebreak
\nl
\mbox{} \vspace{-0.450cm} \\
%

%
% *************************************************************
% *                                                           *
% *              Section 2: Improving K_3, K_4                *
% *                                                           *
% *************************************************************
%

{\bf 2. Derivation of (1.4) and related improvements} \\
% \mbox{\boldmath $K_{\mbox{}_{\!\:\!3}}$}} \\
%
\mbox{} \vspace{-0.750cm} \\

Here we
elaborate
on the method used to obtain (1.3$b$)
in
(\cite{KreissHagstromLorenzZingano2003}, p.\;235)
in order to improve the current estimates
on the regularity time $ \:\!t_{\ast} \!\:\!$
(in dimension $ n = 3, 4$)
as defined in (1.3) above.
We first recall the
elementary Sobolev inequalities \\
\mbox{} \vspace{-0.650cm} \\
\begin{equation}
\tag{2.1$a$}
\|\,\mbox{u}\,
\|_{\mbox{}_{\scriptstyle L^{3}(\mathbb{R}^{3})}}
\leq\:
\Gamma_{\mbox{}_{\!\,\!3}}
\,
\|\,\mbox{u}\,
\|_{\mbox{}_{\scriptstyle L^{2}(\mathbb{R}^{3})}}
  ^{\:\!1/2}
\,\!
\|\, D \:\!\mbox{u}\,
\|_{\mbox{}_{\scriptstyle L^{2}(\mathbb{R}^{3})}}
  ^{\:\!1/2}
\end{equation}
\mbox{} \vspace{-0.650cm} \\
and \\
\mbox{} \vspace{-1.200cm} \\
\begin{equation}
\tag{2.1$b$}
\|\,\mbox{u}\,
\|_{\mbox{}_{\scriptstyle L^{3}(\mathbb{R}^{4})}}
\leq\:
\Gamma_{\mbox{}_{\!\,\!4}}
\,
\|\,\mbox{u}\,
\|_{\mbox{}_{\scriptstyle L^{2}(\mathbb{R}^{4})}}
  ^{\:\!1/3}
\,\!
\|\, D \:\!\mbox{u}\,
\|_{\mbox{}_{\scriptstyle L^{2}(\mathbb{R}^{4})}}
  ^{\:\!2/3}
\end{equation}
\mbox{} \vspace{-0.150cm} \\
for functions in $ H^{1}(\mathbb{R}^{3}) $
and $ H^{1}(\mathbb{R}^{4}) $, respectively.
In both cases,
extremals are given by
$ \;\!\mbox{u}(x) =\;\! C \!\cdot
\mbox{sech}^{2} \:\!(\:\!
\lambda \,|\,x - x_0 \;\!|_{\mbox{}_{2}}) $
for arbitrary $ C \!\in \mathbb{R} $,
$ \lambda \neq 0 $, $ x_0 \!\in \mathbb{R}^{n} \!$
(see \cite{Agueh2008}, p.\;761),
so that the optimal constants in (2.1)
satisfy \\
\mbox{} \vspace{-0.600cm} \\
\begin{equation}
\tag{2.2}
\Gamma_{\mbox{}_{\!\:\!3}} \!\;\!<\,
\mbox{\small $0.558\,901\,115\,737 $},
\qquad
\Gamma_{\mbox{}_{\!\:\!4}} \!\;\!<\,
\mbox{\small $0.419\,577\,519\,172$}.
\end{equation}

\mbox{} \vspace{-0.850cm} \\
We begin with the following lemma. \\
\mbox{} \vspace{-0.050cm} \\
%
% -----------------------------------------------
%
%                 lemma 2.1
%
% -----------------------------------------------
%
{\bf Lemma 2.1.}
\textit{%
Let $\;\! n = 3, 4 $,
and
$ \;\!{\bf u} = (\;\!u_{\mbox{}_{1}} \!\;\!,\!\;\!..., u_{\mbox{}_{\scriptstyle n}})
\in H^{2}(\mathbb{R}^{n})^{\:\!n} \!$.
Then \\
}
\mbox{} \vspace{-0.800cm} \\
\begin{equation}
\tag{2.3$a$}
\int_{\mbox{}_{\scriptstyle \mathbb{R}^{3}}}
\!\!\;\!\Bigl\{\!\!\!\!\:\!
\sum_{\mbox{}\;\;\,i, \,j, \,\ell \,=\,1}^{3}
\!\!\!\!
|\, D_{\ell} \;\!u_{i} \;\!| \:
|\, D_{\ell} \;\!u_{j} \;\!| \:
|\, D_{j} \;\!u_{i} \;\!|
\,\Bigr\}
\;\!dx
\:\leq\:
\Gamma_{\scriptscriptstyle \!\:\!3}^{\;\!3}
\,
\|\, D \:\!{\bf u} \,
\|_{\mbox{}_{\scriptstyle L^{2}(\mathbb{R}^{3})}}
  ^{\:\!3/2}
\|\, D^{2} \,\!{\bf u} \,
\|_{\mbox{}_{\scriptstyle L^{2}(\mathbb{R}^{3})}}
  ^{\:\!3/2}
\end{equation}
\mbox{} \vspace{-0.100cm} \\
\textit{%
if $\;\!n = 3 $,
and \\
}
\mbox{} \vspace{-0.750cm} \\
\begin{equation}
\tag{2.3$b$}
\int_{\mbox{}_{\scriptstyle \mathbb{R}^{4}}}
\!\!\;\!\Bigl\{\!\!\!\!\:\!
\sum_{\mbox{}\;\;\,i, \,j, \,\ell \,=\,1}^{4}
\!\!\!\!
|\, D_{\ell} \;\!u_{i} \;\!| \:
|\, D_{\ell} \;\!u_{j} \;\!| \:
|\, D_{j} \;\!u_{i} \;\!|
\,\Bigr\}
\;\!dx
\:\leq\:
\Gamma_{\scriptscriptstyle \!\:\!4}^{\;\!3}
\,
\|\, D \:\!{\bf u} \,
\|_{\mbox{}_{\scriptstyle L^{2}(\mathbb{R}^{4})}}
\|\, D^{2} \,\!{\bf u} \,
\|_{\mbox{}_{\scriptstyle L^{2}(\mathbb{R}^{4})}}^{\:\!2}
\end{equation}
\mbox{} \vspace{-0.050cm} \\
\textit{%
if $\;\!n = 4 $,
where
$\:\!\Gamma_{\mbox{}_{\!\:\!3}}\!\;\!$,
$\Gamma_{\mbox{}_{\!\:\!4}}\!$
are the constants in
\:\!\mbox{\em (2.1)}, \mbox{\em (2.2)}\!
above.
}

%
% ------------------------------ End of Lemma 2.1
%

%
\mbox{} \vspace{+0.050cm} \\
%
% ----------------------------- Proof of lemma 2.1:
%
{\small
{\bf Proof:}
By repeated application
of the Cauchy-Schwarz inequality,
we have \\
\mbox{} \vspace{-0.100cm} \\
\mbox{} \hfill
$ {\displaystyle
\sum_{\mbox{}\;\;\,i, \,j, \,\ell \,=\,1}^{n}
\!\!\!\!\!\;\!
|\, D_{\ell} \;\!u_{i} \;\!| \:
|\, D_{\ell} \;\!u_{j} \;\!| \:
|\, D_{j} \;\!u_{i} \;\!|
\;\leq\!
\sum_{\mbox{}\;\, i, \,\ell \,=\,1}^{n}
\!\!\:\!
|\, D_{\ell} \;\!u_{i} \;\!|
\:
\Bigl\{\;\!
\sum_{j\,=\,1}^{n}
|\, D_{\ell} \;\!u_{j} \;\!|^{2}
\;\!\Bigr\}^{\!\!\:\!1/2}
\!\;\!
\Bigl\{\,
\sum_{j\,=\,1}^{n}
|\, D_{j} \:\!u_{i} \;\!|^{2}
\;\!\Bigr\}^{\!\!\:\!1/2}
} $ \\
\mbox{} \vspace{+0.050cm} \\
\mbox{} \hfill
$ {\displaystyle
\leq\:
\sum_{i\,=\,1}^{n}
\;
\Bigl\{\,
\sum_{j\,=\,1}^{n}
|\, D_{j} \:\!u_{i} \;\!|^{2}
\;\!\Bigr\}^{\!\!\:\!1/2}
\,\!
\Bigl\{\,
\sum_{\ell\,=\,1}^{n}
|\, D_{\ell} \;\!u_{i} \;\!|^{2}
\;\!\Bigr\}^{\!\!\:\!1/2}
\,\!
\Bigl\{\!
\sum_{\mbox{}\;\,j,\,\ell\,=\,1}^{n}
\!\!\!\;\!
|\, D_{\ell} \;\!u_{j} \;\!|^{2}
\;\!\Bigr\}^{\!\!\:\!1/2}
} $ \\
\mbox{} \vspace{+0.050cm} \\
\mbox{} \hfill
$ {\displaystyle
\leq\;
\Bigl\{\!
\sum_{\mbox{}\;j,\,\ell\,=\,1}^{n}
\!\!\!\;\!
|\, D_{\ell} \;\!u_{j} \;\!|^{2}
\;\!\Bigr\}^{\!\!\:\!1/2}
\:\!
\Bigl\{\!
\sum_{\mbox{}\;i,\,j\,=\,1}^{n}
\!\!\!\;\!
|\, D_{j} \:\!u_{i} \;\!|^{2}
\;\!\Bigr\}^{\!\!\:\!1/2}
\:\!
\Bigl\{\!
\sum_{\mbox{}\;i,\,\ell\,=\,1}^{n}
\!\!\!\;\!
|\, D_{\ell} \;\!u_{i} \;\!|^{2}
\;\!\Bigr\}^{\!\!\:\!1/2}
\!\!\!\!
} $, \\
\mbox{} \vspace{-0.350cm} \\
so that \\
\mbox{} \vspace{-0.700cm} \\
\begin{equation}
\notag
\int_{\mbox{}_{\scriptstyle \mathbb{R}^{n}}}
\!\!\;\!\Bigl\{\!\!\!\!\:\!
\sum_{\mbox{}\;\;\,i, \,j, \,\ell \,=\,1}^{n}
\!\!\!\!
|\, D_{\ell} \;\!u_{i} \;\!| \:
|\, D_{\ell} \;\!u_{j} \;\!| \:
|\, D_{j} \;\!u_{i} \;\!|
\,\Bigr\}
\;\!dx
\:\leq\;
\|\: \mbox{w} \:
\|_{\mbox{}_{\scriptstyle L^{3}(\mathbb{R}^{n})}}^{3}
\!,
\quad \;\,
\mbox{w}(x) :=
\Bigl\{\!\!\;\!
\sum_{\mbox{} \;i, \,j \,=\,1}^{n}
\!\!
|\, D_{j} \;\!u_{i} \;\!|^{\:\!2}
\,
\Bigr\}^{\!\!\;\!1/2}
\!\!\!\!.
\end{equation}
\mbox{} \vspace{+0.025cm} \\
Applying (2.1) to the function
$ \,\mbox{w} \,$
then gives the result,
since,
by (1.5),
we have \\
\mbox{} \vspace{-0.100cm} \\
\mbox{} \hspace{+1.500cm}
$ {\displaystyle
\|\,\mbox{w}\,\|_{\mbox{}_{\scriptstyle L^{2}(\mathbb{R}^{n})}}
=\;
\|\,D\:\!{\bf u}\,\|_{\mbox{}_{\scriptstyle L^{2}(\mathbb{R}^{n})}}
\!\;\!,
\quad
\;\;
\|\,D\:\!\mbox{w}\,\|_{\mbox{}_{\scriptstyle L^{2}(\mathbb{R}^{n})}}
\leq\;\!
\|\,D^{2} \,\!{\bf u}\,\|_{\mbox{}_{\scriptstyle L^{2}(\mathbb{R}^{n})}}
\!\;\!
} $.
\mbox{} \hfill $\Box$ \\
%
%
% -------------------------------- End of Proof of lemma 2.1
%
%
}
\mbox{} \vspace{-0.500cm} \\

We are now in good standing
to reexamine (1.3)
of Section~1.
\!Starting with $ n = 3 $
and
recalling the
Leray's regularity time $ \:\!t_{\ast} $
given in (1.3$a$), (1.3$b$),
we have\:\!: \\
\mbox{} \vspace{-0.010cm} \\
%
% -----------------------------------------------
%
%                 Theorem 2.1
%
% -----------------------------------------------
%
{\bf Theorem 2.1.}
\textit{%
Let $\;\! n = 3 $,
and let $\;\!\mbox{\boldmath $u$}(\cdot,t) $
be any given Leray solution
to the Navier-Stokes system $\,(1.1)$.
Then
there exists $ \;\!t_{\ast\ast} \!\;\!$
satisfying
} \\
\mbox{} \vspace{-0.700cm} \\
\begin{equation}
\tag{2.4}
t_{\ast}
\;\!\leq\,
t_{\ast\ast}
%
% \,\!<\:
%
% % \mbox{\small $ {\displaystyle
% % \frac{1}{\,128\,\pi^{2}} }$}
%
% \mbox{\small $ 0.000\,464\,504\,284 $}
%
\,\!\leq\:
\mbox{\small $ {\displaystyle \frac{1}{2} }$}
\:
\Gamma_{\scriptscriptstyle \!\:\!3}^{\;\!12}
\cdot \;\!
\nu^{-\,5} \,
\|\, \mbox{\boldmath $u$}_0 \;\!
\|_{\mbox{}_{\scriptstyle L^{2}(\mathbb{R}^{3})}}
  ^{\:\!4}
\end{equation}
\mbox{} \vspace{-0.250cm} \\
\textit{%
such that
$ {\displaystyle
\;\!
\|\, D\:\!\mbox{\boldmath $u$}(\cdot,t) \,
\|_{\mbox{}_{\scriptstyle L^{2}(\mathbb{R}^{3})}}
\!\!\:\!
} $
is monotonically decreasing
everywhere in
$ [\,t_{\ast\ast}\!\;\!, \infty) $,
where
$ \Gamma_{\mbox{}_{\!\:\!3}} \!\:\! $
is given in \:\!\mbox{\em (2.1$a$)},
\mbox{\em (2.2)}\!\;\! above.
}
%
% ----------------------------------- End of Theorem 2.1
%

%
\mbox{} \vspace{-0.850cm} \\
%
% ---------------------------------- Proof of Theorem 2.1:
%
{\small
{\bf Proof:}
Consider
$ \:\!\hat{t} > 0 $
satisfying \\
\mbox{} \vspace{-0.625cm} \\
\begin{equation}
\tag{2.5}
\hat{t} \;>\:
\frac{\;\!1\;\!}{2} \:
\Gamma_{\scriptscriptstyle \!\:\!3}^{\;\!12}
\;\!\;\!
\nu^{-\,5}
\,
\|\, \mbox{\boldmath $u$}_{0} \;\!
\|_{\mbox{}_{\scriptstyle L^{2}(\mathbb{R}^{3})}}^{4}
\!\!\!\;\! ~.
\end{equation}
\mbox{} \vspace{-0.350cm} \\
Since
(by (1.2))
$ {\displaystyle
\!\!\;\!
\int_{0}^{\:\!\mbox{\footnotesize $\hat{t}$}}
\!\!\:\!
\|\, D \:\!\mbox{\boldmath $u$}(\cdot,\tau) \,
\|_{\mbox{}_{\scriptstyle L^{2}(\mathbb{R}^{3})}}^{2}
d\tau
\,\!\leq\;\!
\mbox{\small $ {\displaystyle \frac{1}{2\;\!\nu} }$}
\:
\|\, \mbox{\boldmath $u$}_0 \;\!
\|_{\mbox{}_{\scriptstyle L^{2}(\mathbb{R}^{3})}}^{2}
\!\!\;\!
} $,
there exists
some set $ E \!\;\!\subseteq\!\;\! (\;\!0,\!\;\!\;\!\hat{t}\;\!) $
of positive measure
such that \\
\mbox{} \vspace{-0.600cm} \\
\begin{equation}
\tag{2.6}
\|\, D \:\!\mbox{\boldmath $u$}(\cdot,t^{\prime}) \,
\|_{\mbox{}_{\scriptstyle L^{2}(\mathbb{R}^{3})}}
\,\leq\:
\frac{1}{\sqrt{\:\!2\;\!\nu\,}\,}
\:
\|\, \mbox{\boldmath $u$}_{0} \;\!
\|_{\mbox{}_{\scriptstyle L^{2}(\mathbb{R}^{3})}}
\;\!
\hat{t}^{\:-\,1/2}
\quad \;\;\,
\forall \;\,
t^{\:\!\prime} \!\:\!\in E.
\end{equation}
\mbox{} \vspace{-0.150cm} \\
By the {\em epochs of regularity property\/}
\cite{Leray1934} (see also \cite{Heywood1988}),
we can then choose
$ \:\!t^{\:\!\prime} \!\in\!\;\! E $
such that $ \:\!\mbox{\boldmath $u$}(\cdot,\tau) $
is smooth for $\:\!\tau\!\;\!$ close to $ \:\!t^{\:\!\prime} \!$.
Hence,
by (2.5),
we have \\
\mbox{} \vspace{-0.600cm} \\
\begin{equation}
\tag{2.7}
\Gamma_{\mbox{}_{\!\:\!3}}^{\;\!3}
\,
\|\, \mbox{\boldmath $u$}(\cdot,\tau) \,
\|_{\mbox{}_{\scriptstyle L^{2}(\mathbb{R}^{3})}}
  ^{\:\!1/2}
\:\!
\|\, D\:\!\mbox{\boldmath $u$}(\cdot,\tau) \,
\|_{\mbox{}_{\scriptstyle L^{2}(\mathbb{R}^{3})}}
  ^{\:\!1/2}
\!\,\!<\; \nu
\end{equation}
\mbox{} \vspace{-0.175cm} \\
for all $ \tau \geq t^{\:\!\prime} $
close to
the point $ t^{\prime} \!\:\!$.
This gives

%
% differentiating (1.1$a$)
% with respect to $ \:\!x_{\mbox{}_{\scriptstyle \!\;\!\ell}} $,
% multiplying by
% $ D_{\mbox{}_{\scriptstyle \!\:\!\ell}} \mbox{\boldmath $u$}(\cdot,t) $,
% integrating on $\,[\,t^{\prime} \!\:\!, \:\!t\;\!] $,
% and summing over $ \:\!\ell = 1, 2, 3$, \\
%
\mbox{} \vspace{-0.100cm} \\
\mbox{} \hspace{+3.500cm}
$ {\displaystyle
\|\, D \,\!\mbox{\boldmath $u$}(\cdot,t) \,
\|_{\mbox{}_{\scriptstyle L^{2}(\mathbb{R}^{3})}}
  ^{\:\!2}
\!\:\!+\:
2 \,\nu \!\!\;\!
\int_{\mbox{\footnotesize $ t^{\:\!\prime} $}}
    ^{\mbox{\footnotesize $\:\!t$}}
\!
\|\, D^{2} \,\!\mbox{\boldmath $u$}(\cdot,\tau) \,
\|_{\mbox{}_{\scriptstyle L^{2}(\mathbb{R}^{3})}}
  ^{\:\!2}
d\tau
} $ \\
\mbox{} \vspace{-0.625cm} \\
\mbox{} \hfill (2.8) \\
\mbox{} \vspace{-0.425cm} \\
\mbox{} \hfill
$ {\displaystyle
\leq\;\:\!
\|\, D \:\!\mbox{\boldmath $u$}(\cdot,t^{\prime}) \,
\|_{\mbox{}_{\scriptstyle L^{2}(\mathbb{R}^{3})}}
  ^{\:\!2}
\!\:\!+\,
2 \hspace{-0.200cm}
\sum_{i, \, j, \, \ell\,=\,1}^{3}
\int_{\mbox{\footnotesize $ t^{\:\!\prime} $}}
    ^{\mbox{\footnotesize $\:\!t$}}
\!
\int_{\mathbb{R}^{3}}
\!\!\;\!
|\, D_{\mbox{}_{\scriptstyle \!\ell}} \:\! u_{i}(x,\tau) \,|
\;
|\, D_{\mbox{}_{\scriptstyle \!\ell}} \:\! u_{j}(x,\tau) \,|
\;
|\, D_{\scriptstyle \!j} \;\!u_{i}(x,\tau) \,|
\;
dx \, d\tau
} $ \\
\mbox{} \vspace{+0.050cm} \\
\mbox{} \hspace{+0.850cm}
$ {\displaystyle
\leq\;
\|\, D \:\!\mbox{\boldmath $u$}(\cdot,t^{\:\!\prime}) \,
\|_{\mbox{}_{\scriptstyle L^{2}(\mathbb{R}^{3})}}
  ^{\:\!2}
\!\:\!+\,
2 \!\!\;\!
\int_{\mbox{\footnotesize $ t^{\:\!\prime} $}}
    ^{\mbox{\footnotesize $\:\!t$}}
\!
\Gamma_{\mbox{}_{\!\:\!3}}^{\;\!3}
\,
\|\, D \:\!\mbox{\boldmath $u$}(\cdot,\tau) \,
\|_{\mbox{}_{\scriptstyle L^{2}(\mathbb{R}^{3})}}
  ^{\:\!3/2}
\|\, D^{2} \,\!\mbox{\boldmath $u$}(\cdot,\tau) \,
\|_{\mbox{}_{\scriptstyle L^{2}(\mathbb{R}^{3})}}
  ^{\:\!3/2}
\,\!
d\tau
} $
\mbox{} \hfill
\mbox{[}$\,$by (2.3$a$)$\,$\mbox{]} \\
\mbox{} \vspace{+0.025cm} \\
\mbox{} \hfill
$ {\displaystyle
\leq\;
\|\, D \:\!\mbox{\boldmath $u$}(\cdot,t^{\:\!\prime}) \,
\|_{\mbox{}_{\scriptstyle L^{2}(\mathbb{R}^{3})}}
  ^{\:\!2}
\!\:\!+\,
2 \!\!\;\!
\int_{\mbox{\footnotesize $ t^{\:\!\prime} $}}
    ^{\mbox{\footnotesize $\:\!t$}}
\!
\bigl[\;\!\;\!
\Gamma_{\mbox{}_{\!\:\!3}}^{\;\!3}
\,
\|\, \mbox{\boldmath $u$}(\cdot,\tau) \,
\|_{\mbox{}_{\scriptstyle L^{2}(\mathbb{R}^{3})}}
  ^{\:\!1/2}
\|\, D \:\!\mbox{\boldmath $u$}(\cdot,\tau) \,
\|_{\mbox{}_{\scriptstyle L^{2}(\mathbb{R}^{3})}}
  ^{\:\!1/2}
\:\!
\bigr]
\:
\|\, D^{2} \,\!\mbox{\boldmath $u$}(\cdot,\tau) \,
\|_{\mbox{}_{\scriptstyle L^{2}}}
  ^{\:\!2}
d\tau
} $ \\
\mbox{} \vspace{-0.025cm} \\
\mbox{} \hspace{+2.950cm}
$ {\displaystyle
\leq\:
\|\, D \:\!\mbox{\boldmath $u$}(\cdot,t^{\:\!\prime}) \,
\|_{\mbox{}_{\scriptstyle L^{2}(\mathbb{R}^{3})}}
  ^{\:\!2}
\!\:\!+\,
2 \, \nu \!\!\;\!
\int_{\mbox{\footnotesize $ t^{\:\!\prime} $}}
    ^{\mbox{\footnotesize $\:\!t$}}
\!
\|\, D^{2} \,\!\mbox{\boldmath $u$}(\cdot,\tau) \,
\|_{\mbox{}_{\scriptstyle L^{2}(\mathbb{R}^{3})}}
  ^{\:\!2}
d\tau
} $
\mbox{} \hfill
\mbox{[}$\,$by (2.7)$\,$\mbox{]} \\
\mbox{} \vspace{+0.100cm} \\
for all $ \;\!t \geq\,\! t^{\:\!\prime} \!\:\!$
close to $ t^{\:\!\prime} \!\,\!$,
where
in the fourth line above
we used the elementary estimate \\
\mbox{} \vspace{-0.575cm} \\
\begin{equation}
\notag
\|\, D \:\!\mbox{u} \,
\|_{\mbox{}_{\scriptstyle L^{2}(\mathbb{R}^{n})}}
\:\!\leq\;
\|\, \mbox{u} \,
\|_{\mbox{}_{\scriptstyle L^{2}(\mathbb{R}^{n})}}
  ^{\:\!1/2}
\,\!
\|\, D^{2} \mbox{u} \,
\|_{\mbox{}_{\scriptstyle L^{2}(\mathbb{R}^{n})}}
  ^{\:\!1/2}
\end{equation}
\mbox{} \vspace{-0.175cm} \\
(valid for any $\:\!n$),
which is easily obtained
using the Fourier transform.
This shows that \linebreak
$ {\displaystyle
\|\, D \:\!\mbox{\boldmath $u$}(\cdot,t) \,
\|_{\scriptstyle L^{2}(\mathbb{R}^{3})}
\!\;\!
} $
stays bounded by
$ {\displaystyle
\|\, D \:\!\mbox{\boldmath $u$}(\cdot,t^{\:\!\prime}) \,
\|_{\scriptstyle L^{2}(\mathbb{R}^{3})}
\!\;\!
} $,
and because
$ {\displaystyle
\|\, \mbox{\boldmath $u$}(\cdot,t) \,
\|_{\scriptstyle L^{2}(\mathbb{R}^{3})}
\!\;\!
} $
cannot surpass
$ {\displaystyle
\;\!
\|\, \mbox{\boldmath $u$}(\cdot,t^{\:\!\prime}) \,
\|_{\scriptstyle L^{2}(\mathbb{R}^{3})}
\!\;\!
} $
(in view of (1.2)),
it follows
that
we must then have \\
\mbox{} \vspace{-0.550cm} \\
\begin{equation}
\tag{2.9}
\Gamma_{\mbox{}_{\!\:\!3}}^{\;\!3}
\,
\|\, \mbox{\boldmath $u$}(\cdot,t) \,
\|_{\mbox{}_{\scriptstyle L^{2}(\mathbb{R}^{3})}}
  ^{\:\!1/2}
\:\!
\|\, D\:\!\mbox{\boldmath $u$}(\cdot,t) \,
\|_{\mbox{}_{\scriptstyle L^{2}(\mathbb{R}^{3})}}
  ^{\:\!1/2}
\!\;\!<\: \nu
%
% \qquad
% \forall \;\,
% t \geq t^{\:\!\prime}
% \!\:\!.
%
\end{equation}
\mbox{} \vspace{-0.200cm} \\
for {\em all} $ \, t \geq t^{\:\!\prime} \!$,
and in particular
$ \:\!\mbox{\boldmath $u$}(\cdot,t) $
%
% is smooth for $ \:\![\, t^{\:\!\prime} \!, \:\! \infty) $.
%
is smooth for all
$ \;\! t \in [\;\!t^{\:\!\prime} \!\,\!, \infty) $.
So,
all
the estimates
in (2.8) above
can be done
on {\em any\/} interval
$ [\,t_0, \:\!t\, ]  \subseteq
[\,t^{\:\!\prime}\!\:\!, \infty) $,
giving,
for any $ \;\! t \geq t_0 \!\;\!\geq t^{\:\!\prime}$: \\
\mbox{} \vspace{-0.900cm} \\
\mbox{} \hspace{+3.500cm}
$ {\displaystyle
\|\, D \:\!\mbox{\boldmath $u$}(\cdot,t) \,
\|_{\mbox{}_{\scriptstyle L^{2}(\mathbb{R}^{3})}}
  ^{\:\!2}
\!\:\!+\:
2 \, \nu \!\!\;\!
\int_{\mbox{\footnotesize $ t_0 $}}
    ^{\mbox{\footnotesize $\:\!t$}}
\!
\|\, D^{2} \,\!\mbox{\boldmath $u$}(\cdot,\tau) \,
\|_{\mbox{}_{\scriptstyle L^{2}(\mathbb{R}^{3})}}
  ^{\:\!2}
d\tau
} $ \\
\mbox{} \vspace{-0.075cm} \\
\mbox{} \hfill
$ {\displaystyle
\leq\:
\|\, D \:\!\mbox{\boldmath $u$}(\cdot,t_0) \,
\|_{\mbox{}_{\scriptstyle L^{2}(\mathbb{R}^{3})}}
  ^{\:\!2}
\!\:\!+\,
2 \!\!\;\!
\int_{\mbox{\footnotesize $ t_0 $}}
    ^{\mbox{\footnotesize $\:\!t$}}
\!\!\;\!
\bigl[\:
\Gamma_{\mbox{}_{\!\:\!3}}^{\;\!3}
\,
\|\, \mbox{\boldmath $u$}(\cdot,\tau) \,
\|_{\mbox{}_{\scriptstyle L^{2}(\mathbb{R}^{3})}}
  ^{\:\!1/2}
\|\, D \:\!\mbox{\boldmath $u$}(\cdot,\tau) \,
\|_{\mbox{}_{\scriptstyle L^{2}(\mathbb{R}^{3})}}
  ^{\:\!1/2}
\:\!
\bigr]
\:
\|\, D^{2} \,\!\mbox{\boldmath $u$}(\cdot,\tau) \,
\|_{\mbox{}_{\scriptstyle L^{2}(\mathbb{R}^{3})}}
  ^{\:\!2}
\,\!
d\tau
} $ \\
\mbox{} \vspace{-0.075cm} \\
\mbox{} \hspace{+2.950cm}
$ {\displaystyle
\leq\:
\|\, D \:\!\mbox{\boldmath $u$}(\cdot,t_0) \,
\|_{\mbox{}_{\scriptstyle L^{2}(\mathbb{R}^{3})}}
  ^{\:\!2}
\!\:\!+\:
2 \, \nu \!\!\;\!
\int_{\mbox{\footnotesize $ t_0 $}}
    ^{\mbox{\footnotesize $\:\!t$}}
\!\!\;\!
\|\, D^{2} \,\!\mbox{\boldmath $u$}(\cdot,\tau) \,
\|_{\mbox{}_{\scriptstyle L^{2}(\mathbb{R}^{3})}}
  ^{\:\!2}
d\tau
} $.
\mbox{} \hfill
\mbox{[}$\,$by (2.9)$\,$\mbox{]} \\
\mbox{} \vspace{+0.050cm} \\
This shows that
$ {\displaystyle
\;\!
\|\, D\:\! \mbox{\boldmath $u$}(\cdot,t) \,
\|_{\mbox{}_{\scriptstyle L^{2}(\mathbb{R}^{3})}}
\!
} $
stays finite and is monotonically decreasing
% for $ \:\!t \geq t^{\:\!\prime} $
everywhere in
$ {\displaystyle
[\,t^{\:\!\prime}\!, \:\!\infty)
\supseteq\:\!
[\,\hat{t}, \:\!\infty)
} $.
$\!$Therefore,
by Leray's theory,
(1.3$a$) of Section~1
will be satisfied for
any
$ \;\!t_{\ast} \!> t^{\:\!\prime} \!\;\!$,
with
$ \;\!t^{\:\!\prime} \!\:\!<\:\! \hat{t\:\!} \!\;\!$.
Recalling (2.5),
this completes the proof of
{\sc Theorem 2.1}.
%
% since
% $ {\displaystyle
% \:\!
% 1/2
% \;\:\!
% \Gamma_{\mbox{}_{\!\:\!3}}^{\;\!12}
% \!\;\!<\:\!
% 0.000\,464\,504\,284
% } $.
%
%
% ------------------------- End of Proof of Theorem 2.1
%
}
\hfill \mbox{\small $\Box$} \\
%
%
% --------------------------------------------- Theorem 2.2:
%
%
\nl
In much the same way,
using (2.3$b$)
and the 4\mbox{\small D}
version of (2.8),
we can show$\:\!$: \\
%
% we obtain
% the following result. \\
%
%
\mbox{} \vspace{-0.050cm} \\
%
% -----------------------------------------------
%
%                 Theorem 2.2
%
% -----------------------------------------------
%
{\bf Theorem 2.2.}
\textit{%
Let $\;\! n = 4 $,
and let $\;\!\mbox{\boldmath $u$}(\cdot,t) $
be any given Leray solution
to the Navier-Stokes system $\,(1.1)$.
Let $\;\!t_{\ast} \!\;\!\geq 0 $
be defined in $\;\!(1.3$a$)$.
Then
there exists $ \;\!t_{\ast\ast} \!\;\!$
satisfying
} \\
\mbox{} \vspace{-0.650cm} \\
\begin{equation}
\tag{2.10}
t_{\ast}
\;\!\leq\,
t_{\ast\ast}
%
% \,\!<\:
%
% \mbox{\small $ 0.002\,727\,993\,110 $}
%
\:\!\leq\:
\mbox{\small $ {\displaystyle \frac{1}{2} }$}
\:
\Gamma_{\scriptscriptstyle \!\:\!4}^{\;\!6}
\cdot \;\!
\nu^{-\,3}
\,
\|\, \mbox{\boldmath $u$}_0 \;\!
\|_{\mbox{}_{\scriptstyle L^{2}(\mathbb{R}^{4})}}
  ^{\:\!2}
\end{equation}
\mbox{} \vspace{-0.275cm} \\
\textit{%
such that
$ {\displaystyle
\;\!
\|\, D\:\!\mbox{\boldmath $u$}(\cdot,t) \,
\|_{\mbox{}_{\scriptstyle L^{2}(\mathbb{R}^{4})}}
\!\!\:\!
} $
is monotonically decreasing
everywhere in
$ [\,t_{\ast\ast}\!\;\!, \infty) $,
where
$ \Gamma_{\mbox{}_{\!\:\!4}} \!\,\! $
is given in \,\!\mbox{\em (2.1$b$)},
\!\:\!\mbox{\em (2.2)}\!\:\! above.
}
%
% ----------------------------------- End of Theorem 2.2
%

%
\mbox{} \vspace{-0.050cm} \\
%
% --------------------------------- Proof of Theorem 2.2:
%
{\small
{\bf Proof:}
Consider
$ \:\!\hat{t} > 0 $
satisfying \\
\mbox{} \vspace{-0.625cm} \\
\begin{equation}
\tag{2.11}
\hat{t} \;>\:
\frac{\;\!1\;\!}{2} \:
\Gamma_{\scriptscriptstyle \!\:\!4}^{\;\!6}
\;\!\;\!
\nu^{-\,3}
\,
\|\, \mbox{\boldmath $u$}_{0} \;\!
\|_{\mbox{}_{\scriptstyle L^{2}(\mathbb{R}^{4})}}^{2}
\!\!\!\;\! ~.
\end{equation}
\mbox{} \vspace{-0.350cm} \\
As before,
from (1.2)
it follows the existence
of some set
$ E \!\;\!\subseteq\!\;\! (\;\!0,\!\;\!\;\!\hat{t}\;\!) $
with positive measure
such that \\
\mbox{} \vspace{-0.800cm} \\
\begin{equation}
\tag{2.12}
\|\, D \:\!\mbox{\boldmath $u$}(\cdot,t^{\prime}) \,
\|_{\mbox{}_{\scriptstyle L^{2}(\mathbb{R}^{4})}}
\,\leq\:
\frac{1}{\sqrt{\:\!2\;\!\nu\,}\,}
\:
\|\, \mbox{\boldmath $u$}_{0} \;\!
\|_{\mbox{}_{\scriptstyle L^{2}(\mathbb{R}^{4})}}
\;\!
\hat{t}^{\:-\,1/2}
\quad \;\;\,
\forall \;\,
t^{\:\!\prime} \!\:\!\in E
\end{equation}
\mbox{} \vspace{-0.150cm} \\
and so,
by the epochs of regularity property,
we can again choose
$ \;\!t^{\:\!\prime} \!\in\!\;\! E \;\!$
such that $ \:\!\mbox{\boldmath $u$}(\cdot,\tau) $
is smooth for $\:\!\tau\!\;\!$ close to $ \:\!t^{\:\!\prime} \!$.
Hence,
by (2.11) and (2.12),
we have \\
\mbox{} \vspace{-0.600cm} \\
\begin{equation}
\tag{2.13}
\Gamma_{\mbox{}_{\!\:\!4}}^{\;\!3}
\;\!\;\!
\|\, D\:\!\mbox{\boldmath $u$}(\cdot,\tau) \,
\|_{\mbox{}_{\scriptstyle L^{2}(\mathbb{R}^{4})}}
\!\;\!<\; \nu
\end{equation}
\mbox{} \vspace{-0.175cm} \\
for all $ \tau \geq t^{\:\!\prime} $
close to
the point $ t^{\prime} \!\:\!$.
This then gives

\mbox{} \vspace{-0.100cm} \\
\mbox{} \hspace{+3.500cm}
$ {\displaystyle
\|\, D \,\!\mbox{\boldmath $u$}(\cdot,t) \,
\|_{\mbox{}_{\scriptstyle L^{2}(\mathbb{R}^{4})}}
  ^{\:\!2}
\!\:\!+\:
2 \,\nu \!\!\;\!
\int_{\mbox{\footnotesize $ t^{\:\!\prime} $}}
    ^{\mbox{\footnotesize $\:\!t$}}
\!
\|\, D^{2} \,\!\mbox{\boldmath $u$}(\cdot,\tau) \,
\|_{\mbox{}_{\scriptstyle L^{2}(\mathbb{R}^{4})}}
  ^{\:\!2}
d\tau
} $ \\
\mbox{} \vspace{-0.025cm} \\
\mbox{} \hfill
$ {\displaystyle
\leq\;\:\!
\|\, D \:\!\mbox{\boldmath $u$}(\cdot,t^{\prime}) \,
\|_{\mbox{}_{\scriptstyle L^{2}(\mathbb{R}^{4})}}
  ^{\:\!2}
\!\:\!+\,
2 \hspace{-0.200cm}
\sum_{i, \, j, \, \ell\,=\,1}^{4}
\int_{\mbox{\footnotesize $ t^{\:\!\prime} $}}
    ^{\mbox{\footnotesize $\:\!t$}}
\!
\int_{\mathbb{R}^{4}}
\!\!\;\!
|\, D_{\mbox{}_{\scriptstyle \!\ell}} \:\! u_{i}(x,\tau) \,|
\;
|\, D_{\mbox{}_{\scriptstyle \!\ell}} \:\! u_{j}(x,\tau) \,|
\;
|\, D_{\scriptstyle \!j} \;\!u_{i}(x,\tau) \,|
\;
dx \, d\tau
} $

\newpage
\mbox{} \vspace{-0.850cm} \\
\mbox{} \hspace{+0.850cm}
$ {\displaystyle
\leq\;
\|\, D \:\!\mbox{\boldmath $u$}(\cdot,t^{\:\!\prime}) \,
\|_{\mbox{}_{\scriptstyle L^{2}(\mathbb{R}^{4})}}
  ^{\:\!2}
\!\:\!+\,
2 \!\!\;\!
\int_{\mbox{\footnotesize $ t^{\:\!\prime} $}}
    ^{\mbox{\footnotesize $\:\!t$}}
\!
\Gamma_{\mbox{}_{\!\:\!4}}^{\;\!3}
\,
\|\, D \:\!\mbox{\boldmath $u$}(\cdot,\tau) \,
\|_{\mbox{}_{\scriptstyle L^{2}(\mathbb{R}^{4})}}
\|\, D^{2} \,\!\mbox{\boldmath $u$}(\cdot,\tau) \,
\|_{\mbox{}_{\scriptstyle L^{2}(\mathbb{R}^{4})}}
  ^{\:\!2}
\,\!
d\tau
} $
\mbox{} \hfill
\mbox{[}$\,$by (2.3$b$)$\,$\mbox{]} \\
\mbox{} \vspace{-0.025cm} \\
\mbox{} \hspace{+2.950cm}
$ {\displaystyle
\leq\:
\|\, D \:\!\mbox{\boldmath $u$}(\cdot,t^{\:\!\prime}) \,
\|_{\mbox{}_{\scriptstyle L^{2}(\mathbb{R}^{4})}}
  ^{\:\!2}
\!\:\!+\,
2 \, \nu \!\!\;\!
\int_{\mbox{\footnotesize $ t^{\:\!\prime} $}}
    ^{\mbox{\footnotesize $\:\!t$}}
\!
\|\, D^{2} \,\!\mbox{\boldmath $u$}(\cdot,\tau) \,
\|_{\mbox{}_{\scriptstyle L^{2}(\mathbb{R}^{4})}}
  ^{\:\!2}
d\tau
} $
\mbox{} \hfill
\mbox{[}$\,$by (2.13)$\,$\mbox{]} \\
\mbox{} \vspace{-0.350cm} \\
\mbox{} \hfill (2.14) \\
\mbox{} \vspace{-0.375cm} \\
for all $ \;\!t \geq\,\! t^{\:\!\prime} \!\;\!$
close to $ t^{\:\!\prime} \!\,\!$.
This actually implies
that
$ {\displaystyle
\;\!
\|\, D \:\!\mbox{\boldmath $u$}(\cdot,t) \,
\|_{\scriptstyle L^{2}(\mathbb{R}^{4})}
\!\;\!
} $
will stay bounded by
$ {\displaystyle
\|\, D \:\!\mbox{\boldmath $u$}(\cdot,t^{\:\!\prime}) \,
\|_{\scriptstyle L^{2}(\mathbb{R}^{4})}
\!\;\!
} $
for {\em all}
$ \,t \geq\,\! t^{\:\!\prime} \!\,\!$,
\;\!so that
we actually have \\
\mbox{} \vspace{-0.550cm} \\
\begin{equation}
\tag{2.15}
\Gamma_{\mbox{}_{\!\:\!4}}^{\;\!3}
\;\!\;\!
\|\, D\:\!\mbox{\boldmath $u$}(\cdot,\tau) \,
\|_{\mbox{}_{\scriptstyle L^{2}(\mathbb{R}^{4})}}
\!\;\!<\: \nu
\qquad
\forall \;\,
\tau \geq t^{\:\!\prime}
\!\:\!.
\end{equation}
\mbox{} \vspace{-0.200cm} \\
In particular,
it follows that
$ \;\!\mbox{\boldmath $u$}(\cdot,t) \:\!$
is smooth for every
$ \;\! t \in \:\![\, t^{\:\!\prime} \!, \:\! \infty) $.
This fact
and (2.15)
allow us to repeat
the derivation of
(2.14) above
on {\em any\/} interval
$ [\,t_0, \:\!t\, ]  \subseteq
[\,t^{\:\!\prime}\!\:\!, \infty) $,
obtaining \\
\mbox{} \vspace{-0.100cm} \\
\mbox{} \hspace{+3.500cm}
$ {\displaystyle
\|\, D \:\!\mbox{\boldmath $u$}(\cdot,t) \,
\|_{\mbox{}_{\scriptstyle L^{2}(\mathbb{R}^{4})}}
  ^{\:\!2}
\!\:\!+\:
2 \, \nu \!\!\;\!
\int_{\mbox{\footnotesize $ t_0 $}}
    ^{\mbox{\footnotesize $\:\!t$}}
\!
\|\, D^{2} \,\!\mbox{\boldmath $u$}(\cdot,\tau) \,
\|_{\mbox{}_{\scriptstyle L^{2}(\mathbb{R}^{4})}}
  ^{\:\!2}
d\tau
} $ \\
\mbox{} \vspace{-0.075cm} \\
\mbox{} \hspace{+1.500cm}
$ {\displaystyle
\leq\:
\|\, D \:\!\mbox{\boldmath $u$}(\cdot,t_0) \,
\|_{\mbox{}_{\scriptstyle L^{2}(\mathbb{R}^{4})}}
  ^{\:\!2}
\!\:\!+\,
2 \!\!\;\!
\int_{\mbox{\footnotesize $ t_0 $}}
    ^{\mbox{\footnotesize $\:\!t$}}
\!\!\;\!
\Gamma_{\mbox{}_{\!\:\!4}}^{\;\!3}
\;\!\;\!
\|\, D \:\!\mbox{\boldmath $u$}(\cdot,\tau) \,
\|_{\mbox{}_{\scriptstyle L^{2}(\mathbb{R}^{4})}}
\;\!
\|\, D^{2} \,\!\mbox{\boldmath $u$}(\cdot,\tau) \,
\|_{\mbox{}_{\scriptstyle L^{2}(\mathbb{R}^{4})}}
  ^{\:\!2}
\,\!
d\tau
} $ \\
\mbox{} \vspace{-0.075cm} \\
\mbox{} \hspace{+2.975cm}
$ {\displaystyle
\leq\:
\|\, D \:\!\mbox{\boldmath $u$}(\cdot,t_0) \,
\|_{\mbox{}_{\scriptstyle L^{2}(\mathbb{R}^{4})}}
  ^{\:\!2}
\!\:\!+\:
2 \, \nu \!\!\;\!
\int_{\mbox{\footnotesize $ t_0 $}}
    ^{\mbox{\footnotesize $\:\!t$}}
\!\!\;\!
\|\, D^{2} \,\!\mbox{\boldmath $u$}(\cdot,\tau) \,
\|_{\mbox{}_{\scriptstyle L^{2}(\mathbb{R}^{4})}}
  ^{\:\!2}
d\tau
} $.
\mbox{} \hfill
\mbox{[}$\,$by (2.15)$\,$\mbox{]} \\
\mbox{} \vspace{+0.050cm} \\
This shows that
$ {\displaystyle
\;\!
\|\, D\:\! \mbox{\boldmath $u$}(\cdot,t) \,
\|_{\mbox{}_{\scriptstyle L^{2}(\mathbb{R}^{4})}}
\!
} $
stays bounded and monotonically decreasing
% for $ \:\!t \geq t^{\:\!\prime} $
everywhere in
$ {\displaystyle
[\,t^{\:\!\prime}\!, \:\!\infty)
\supseteq\:\!
[\,\hat{t}, \:\!\infty)
} $.
$\!$As in the previous proof,
(1.3$a$)
must then be satisfied for
any
$ t_{\ast} \!> t^{\:\!\prime} \!\;\!$,
where
$ \;\!t^{\:\!\prime} \!\:\!<\:\! \hat{t\:\!} \!\;\!$.
Recalling (2.11),
this completes the proof of
{\sc Theorem 2.2}.
%
%
% ------------------------- End of Proof of Theorem 2.2
%
}
\hfill \mbox{\small $\Box$} \\
\nl
%
% -------------------------------------------------------------
%
%                          Remark 2.1
%
% -------------------------------------------------------------
%
{\bf Remark 2.1.}
\!\:\!In dimension $ \:\! n = 2 $,
it is well known that
$ \:\!\mbox{\boldmath $u$} \in C^{\infty}(\,\!\mathbb{R}^{2}
\!\times\!\:\! (\:\!0, \infty)) \,\!$
and \linebreak
\mbox{} \vspace{-0.625cm} \\
\begin{equation}
\tag{2.16}
\mbox{\boldmath $u$}(\cdot,t)
\in
C^{0}(\,\!(\:\!0, \infty), \:\!
\mbox{\boldmath $H$}\mbox{}^{\!\:\!m}\!\;\!(\mathbb{R}^{2}))
\quad \;\,
\forall \;\;\! m \geq 0,
\end{equation}
\mbox{} \vspace{-0.600cm} \\
with \\
\mbox{} \vspace{-0.900cm} \\
\begin{equation}
\tag{2.17}
\|\, D \:\!\mbox{\boldmath $u$}(\cdot,t) \,
\|_{\mbox{}_{\scriptstyle L^{2}(\mathbb{R}^{2})}}
  ^{\:\!2}
\!\:\!+\:
2 \, \nu \!\!\;\!
\int_{\mbox{\footnotesize $ t_0 $}}
    ^{\mbox{\footnotesize $\:\!t$}}
\!
\|\, D^{2} \,\!\mbox{\boldmath $u$}(\cdot,\tau) \,
\|_{\mbox{}_{\scriptstyle L^{2}(\mathbb{R}^{2})}}
  ^{\:\!2}
d\tau
\;=\;
\|\, D \:\!\mbox{\boldmath $u$}(\cdot,t_0) \,
\|_{\mbox{}_{\scriptstyle L^{2}(\mathbb{R}^{2})}}
  ^{\:\!2}
\end{equation}
\mbox{} \vspace{-0.125cm} \\
for all
$ \;\! t \geq t_0 \!\:\!> 0 $,
so that,
in particular,
$ {\displaystyle
\;\!
\|\, D \:\!\mbox{\boldmath $u$}(\cdot,t) \,
\|_{\mbox{}_{\scriptstyle L^{2}(\mathbb{R}^{2})}}
\!\;\!
} $
is monotonically decreasing
in $(\:\!0, \infty) $.
Since
monotonic functions
$ f \!\;\!\in L^{1}(\:\!t_0, \infty) $
satisfy:
$ f(t) = o(1/t) $
as $ t \rightarrow \infty $
(because, in this case:
$ \;\!t \,| \:\!f(t)\;\!| \:\!\leq\,
\mbox{\small 2} \!\:\!\int_{\:\!t/2}^{\;\!t}
| \:\!f(\tau)\;\! | \, d\tau \,\!$
for all $ \:\!t \geq \mbox{\small 2} \;\!t_{0} $),
it follows
that \\
%
% $ {\displaystyle
% \lim_{t\,\rightarrow\,\infty}
% t^{1/2} \,
% \|\, D \:\!\mbox{\boldmath $u$}(\cdot,t) \,
% \|_{\mbox{}_{\scriptstyle L^{2}(\mathbb{R}^{n})}}
% \!=\;\! 0
% } $
% ($ 2 \leq n \leq 4 $),
%
\mbox{} \vspace{-0.600cm} \\
\begin{equation}
\tag{2.18}
\lim_{t\,\rightarrow\,\infty}
\,
t^{1/2} \,
\|\, D \:\!\mbox{\boldmath $u$}(\cdot,t) \,
\|_{\mbox{}_{\scriptstyle L^{2}(\mathbb{R}^{n})}}
\;\!=\; 0
\end{equation}
\mbox{} \vspace{-0.150cm} \\
if $ \;\!2 \leq n \leq 4 $,
as partially observed in
\cite{KreissHagstromLorenzZingano2002, %
KreissHagstromLorenzZingano2003}.
\!(\:\!For a derivation of (2.18)
in the case of more general $n$,
see \cite{OliverTiti2000, SchonbekWiegner1996}.)
\newpage
\mbox{} \vspace{-0.750cm} \\
%
% -------------------------------------------------------------
%
%                          Remark 2.2
%
% -------------------------------------------------------------
%
{\bf Remark 2.2.}
From (2.18),
one can easily obtain
that
all Leray solutions to (1.1)
satisfy the asymptotic property \\
\mbox{} \vspace{-0.650cm} \\
\begin{equation}
\tag{2.19}
\lim_{t\,\rightarrow\,\infty}
\,
\|\, \mbox{\boldmath $u$}(\cdot,t) \,
\|_{\mbox{}_{\scriptstyle L^{2}(\mathbb{R}^{n})}}
\;\!=\; 0
\end{equation}
\mbox{} \vspace{-0.125cm} \\
if $ \;\!2 \leq n \leq 4 $,
as shown in
\cite{RigeloSchutzZinganoZingano2016, Zingano2015}.\footnote{%
%
% -------------------------------------- Footnote 2:
%
The validity or not of (2.19)
was left open by Leray in \cite{Leray1934} (p.\;248),
being first shown in \cite{Kato1984, Masuda1984}.
}
%
% ----------------------------------- End of Footnote 2
%
In fact,
in view that
$ \mbox{\boldmath $u$}(\cdot,t) $
is smooth
for large~\mbox{$ \:\!t $},
it can be written as \\
\mbox{} \vspace{-0.700cm} \\
\begin{equation}
\tag{2.20}
\mbox{} \hspace{+0.500cm}
\mbox{\boldmath $u$}(\cdot,t)
\;=\;
e^{\:\!\nu\:\!\Delta (t - t_0)}
\:\!
\mbox{\boldmath $u$}(\cdot,t_0)
\,-
\int_{\!\:\!t_0}^{\;\!t}
\!\!\:\!
e^{\:\!\nu\:\!\Delta (t - \tau)}
\:\!
\mbox{\boldmath $Q$}(\cdot,\tau)
\,
d\tau,
\quad \;\;\,
t > t_0
\end{equation}
\mbox{} \vspace{-0.150cm} \\
for $ \:\!t_0 \!\;\!$
large enough,
where
$ {\displaystyle
\:\!
\mbox{\boldmath $Q$} \!\;\!=\,\!
\mbox{\boldmath $u$} \!\:\!\cdot\!\:\!\nabla \mbox{\boldmath $u$}
\:\!+\!\;\! \nabla p
\;\!
} $
is the Helmholtz projection
of
$ \,\!\,\!\mbox{\boldmath $u$} \!\;\!\cdot\!\;\!\nabla \mbox{\boldmath $u$} $
in
$ \mbox{\boldmath $L$}^{2}_{\sigma}(\mathbb{R}^{n}) $,
and
$ \:\!e^{\:\!\nu\:\!\Delta t} \!\;\!$
denotes the heat semigroup.
%
%
% -------------------------- Proof of (2.19) for n = 4:
%
%
\!\:\!In the case $ n = 4 $,
we can then
get~(2.19) as follows:
from (2.20),
we obtain,
recalling
that
$ {\displaystyle
\,
\|\, \mbox{u} \,
\|_{\mbox{}_{\scriptstyle L^{4}(\mathbb{R}^{4})}}
\!\;\!\leq\;\!
\|\, D \:\!\mbox{u} \,
\|_{\mbox{}_{\scriptstyle L^{2}(\mathbb{R}^{4})}}
} $ \linebreak
\mbox{} \vspace{-0.600cm} \\
for arbitrary
$ \;\!\mbox{u} \in H^{1}(\mathbb{R}^{4}) $, \\
\mbox{} \vspace{-0.100cm} \\
\mbox{} \hspace{+0.250cm}
$ {\displaystyle
\|\, \mbox{\boldmath $u$}(\cdot,t) \,
\|_{\mbox{}_{\scriptstyle L^{2}(\mathbb{R}^{4})}}
\leq\;
\|\, \mbox{\boldmath $v$}_{0}(\cdot,t) \,
\|_{\mbox{}_{\scriptstyle L^{2}(\mathbb{R}^{4})}}
\:\!+
\int_{\!\:\!t_0}^{\;\!t}
\!\!\:\!
\|\, \mbox{\boldmath $u$}(\cdot,\tau) \!\;\!\cdot\!\;\!
\nabla \mbox{\boldmath $u$}(\cdot,\tau) \,
\|_{\mbox{}_{\scriptstyle L^{2}(\mathbb{R}^{4})}}
\;\!
d\tau
} $ \\
\mbox{} \vspace{-0.000cm} \\
\mbox{} \hspace{+3.025cm}
$ {\displaystyle
\leq\;
\|\, \mbox{\boldmath $v$}_{0}(\cdot,t) \,
\|_{\mbox{}_{\scriptstyle L^{2}(\mathbb{R}^{4})}}
\:\!+\,
\sqrt{\:\!2\;\!\:\!} \!\!
\int_{\!\:\!t_0}^{\;\!t}
\!\!\:\!
\|\, \mbox{\boldmath $u$}(\cdot,\tau) \,
\|_{\mbox{}_{\scriptstyle L^{4}(\mathbb{R}^{4})}}
\,\!
\|\, D \:\! \mbox{\boldmath $u$}(\cdot,\tau) \,
\|_{\mbox{}_{\scriptstyle L^{4}(\mathbb{R}^{4})}}
\;\!
d\tau
} $ \\
\mbox{} \vspace{-0.000cm} \\
\mbox{} \hspace{+3.025cm}
$ {\displaystyle
\leq\;
\|\, \mbox{\boldmath $v$}_{0}(\cdot,t) \,
\|_{\mbox{}_{\scriptstyle L^{2}(\mathbb{R}^{4})}}
\:\!+\,
\sqrt{\:\!2\;\!\:\!} \!\!
\int_{\!\:\!t_0}^{\;\!t}
\!\!\:\!
\|\, D \:\!\mbox{\boldmath $u$}(\cdot,\tau) \,
\|_{\mbox{}_{\scriptstyle L^{2}(\mathbb{R}^{4})}}
\,\!
\|\, D^{2} \,\! \mbox{\boldmath $u$}(\cdot,\tau) \,
\|_{\mbox{}_{\scriptstyle L^{2}(\mathbb{R}^{4})}}
\;\!
d\tau
} $ \\
\mbox{} \vspace{+0.150cm} \\
since
$ {\displaystyle
\;\!
\|\, \mbox{\boldmath $Q$}(\cdot,\tau) \,
\|_{\scriptstyle L^{2}(\mathbb{R}^{n})}
\leq\;\!
\|\, \mbox{\boldmath $u$}(\cdot,\tau) \!\:\!\cdot\!\:\!
\nabla \mbox{\boldmath $u$}(\cdot,\tau) \,
\|_{\scriptstyle L^{2}(\mathbb{R}^{n})}
} $
by the orthogonality of the
Helm\-holtz projector in
$ \mbox{\boldmath $L$}^{2}(\mathbb{R}^{n}) $
(see e.g.\;\cite{RigeloSchutzZinganoZingano2016}),
or directly using Fourier transform
\cite{KreissHagstromLorenzZingano2002, %
KreissHagstromLorenzZingano2003},
and where
$ {\displaystyle
\:\!
\mbox{\boldmath $v$}_{0}(\cdot,t)
\!:=
e^{\:\!\nu\:\!\Delta (t - t_0)}
\:\!
\mbox{\boldmath $u$}(\cdot,t_0)
} $.
\!This shows that,
given $ \epsilon > 0 $,
taking $ \:\!t_0 $
large enough
we get
$ {\displaystyle
\;\!
\|\, \mbox{\boldmath $u$}(\cdot,t) \,
\|_{\scriptstyle L^{2}(\mathbb{R}^{4})}
\!\;\!< \epsilon
\;\!
} $
for all $ t > t_{0}$,
since
the integrand on the
righthand side above is in
$ L^{1}(\:\!t_{\ast}, \infty) $.
A similar argument for
$ n = 2, \,3 $
can be found in
\cite{RigeloSchutzZinganoZingano2016, %
Zingano2015}.
The proof of (2.19)
for arbitrary $ \;\!n \;\!$
is significantly harder
and given in \cite{Wiegner1987}.
\mbox{} \hfill $ \Box $ \\
\nl
\mbox{} \vspace{-0.450cm} \\
%
% -------------------------------------------------------------
%
%                          Remark 2.3
%
% -------------------------------------------------------------
%
{\bf Remark 2.3.}
\!\;\!Now that (2.18) and (2.19) are known,
it is possible to extend~the
monotonicity property
of
$ {\displaystyle
\;\!
\|\, D\,\!\mbox{\boldmath $u$}(\cdot,t) \,
\|_{L^{2}(\mathbb{R}^{n})}
} $,
$ n \leq 4 $
(cf.\;(2.17) and
\mbox{\small \sc Theorems 2.1},~\mbox{\small \sc 2.2})
to higher order derivatives.
Let $ \;\!t_{\ast} $ be the regularity time
given in (1.3) for $ \;\!n = 3, \;\!4 $,
and $ \;\! t_{\ast} \!:= 0 \:\!$ if $ \:\!n = 2 $
(cf.\;(2.16) above).
We then have\:\!: \\
\mbox{} \vspace{+0.025cm} \\
%
% -----------------------------------------------
%
%                 Theorem 2.3
%
% -----------------------------------------------
%
{\bf Theorem 2.3.}
\textit{%
Let $\;\! 2 \leq n \leq 4 $,
and let $\;\!\mbox{\boldmath $u$}(\cdot,t) $
be any particular Leray solution
to the Navier-Stokes system $\,(1.1)$.
Then,
for each $ \;\! m \geq 1 \!\!\:\!:$
there exists
$ \;\!\:\!t_{\ast\ast}^{\:\!(m)} \!\geq\:\! t_{\ast} \!\;\!$
such that
$ {\displaystyle
\,
\|\, D^{m} \mbox{\boldmath $u$}(\cdot,t) \,
\|_{\mbox{}_{\scriptstyle L^{2}(\mathbb{R}^{n})}}
\!\!\:\!
} $
is monotonically decreasing
everywhere in
$ \;\![\,t_{\ast\ast}^{\:\!(m)}\!, \:\!\infty) $.
} \\
%
% ---------------------------------------- Proof of Theorem 2.3:
%
\nl
{\small
{\bf Proof\/:}
Let $\;\! t_0 \!\;\!> t_{\ast} $,
$ m \geq 1 $.
Given $ \;\! t > t_0 $,
we have from (1.1$a$), (1.1$b$)
the energy estimate \\
\mbox{} \vspace{-0.050cm} \\
\mbox{} \hspace{+0.750cm}
$ {\displaystyle
\|\, D^{m} \,\!\mbox{\boldmath $u$}(\cdot,t) \,
\|_{\mbox{}_{\scriptstyle L^{2}(\mathbb{R}^{n})}}
  ^{\:\!2}
\!\:\!+\;
2 \, \nu \!\!\;\!
\int_{\mbox{}_{\scriptstyle \!\:\!t_{\mbox{}_{0}}}}^{\;\!t}
\!\!\:\!
\|\, D^{m + 1} \,\!\mbox{\boldmath $u$}(\cdot,\tau) \,
\|_{\mbox{}_{\scriptstyle L^{2}(\mathbb{R}^{n})}}
  ^{\:\!2}
d\tau
\;\leq\;\:\!
\|\, D^{m} \,\!\mbox{\boldmath $u$}(\cdot,t_0) \,
\|_{\mbox{}_{\scriptstyle L^{2}(\mathbb{R}^{n})}}
  ^{\:\!2}
} $ \\
\mbox{} \vspace{-0.550cm} \\
\mbox{} \hfill (2.21) \\
\mbox{} \vspace{-0.550cm} \\
\mbox{} \hspace{+1.250cm}
$ {\displaystyle
+\;\,
2 \,
% \sum_{\mbox{} \;\;i, \,j, \,j_{\mbox{}_{1}} \!,..., \,j_{\mbox{}_{m}} =\,1}^{n}
\sum
\int_{\mbox{}_{\scriptstyle \!\:\!t_{\mbox{}_{0}}}}^{\;\!t}
\!\!\;\!
\int_{\mbox{}_{\scriptstyle \:\!\mathbb{R}^{n}}}
\!\!\!\:\!
|\, D_{\!\;\!j} \;\!
D_{\!\;\!j_{\mbox{}_{1}}}
\!\!\!\;\!\cdot\!\,\!\cdot\!\,\!\cdot \!\:\!
D_{\!\;\!j_{\mbox{}_{m}}}
\!\:\! u_{i}(x,\tau)\,|
\cdot
|\, D_{\!\;\!j_{\mbox{}_{1}}}
\!\!\!\;\!\cdot\!\,\!\cdot\!\,\!\cdot \!\:\!
D_{\!\;\!j_{\mbox{}_{m}}}
\!\:\! (\:\! u_{i}(x,\tau)\, u_{j}(x,\tau) \,\!)\,|
\;dx \,d\tau
} $ \\
\mbox{} \vspace{+0.075cm} \\
where the sum
is over
all indices
$ \;\!1 \leq\:\! i, \;\!j, \;\!j_{\mbox{}_{1}}\!\:\!,
 ...\;\!, \;\!j_{\mbox{}_{m}} \!\:\!\leq\:\! n $.
Another important inequality is \\
\mbox{} \vspace{-0.550cm} \\
\begin{equation}
\tag{2.22}
\|\, D^{\ell} \,\!\mbox{u} \,
\|_{\mbox{}_{\scriptstyle L^{2}(\mathbb{R}^{n})}}
\:\!\leq\;
\|\, \mbox{u} \,
\|_{\mbox{}_{\scriptstyle L^{2}(\mathbb{R}^{n})}}
  ^{\:\!1 - \theta}
\,\!
\|\, D^{m} \,\!\mbox{u} \,
\|_{\mbox{}_{\scriptstyle L^{2}(\mathbb{R}^{n})}}
  ^{\:\!\theta}
\!\:\!,
\quad \;\;
\theta \,=\:
\mbox{\small $ {\displaystyle \frac{\ell}{m} }$}
\end{equation}
\mbox{} \vspace{-0.100cm} \\
for every $\;\!0 \leq \ell \leq m $,
which can be obtained by Fourier transform.
From this point further, \linebreak
the argument becomes dependent on the dimension $\:\!n$,
and we illustrate the typical steps
by considering, say,
$ \:\! n = 4 $.
In this case,
we observe that,
from (2.21) and H\"older's inequality, \\
\mbox{} \vspace{-0.350cm} \\
\mbox{} \hspace{+0.750cm}
$ {\displaystyle
\|\, D^{m} \,\!\mbox{\boldmath $u$}(\cdot,t) \,
\|_{\mbox{}_{\scriptstyle L^{2}(\mathbb{R}^{4})}}
  ^{\:\!2}
\!\:\!+\;
2 \, \nu \!\!\;\!
\int_{\mbox{}_{\scriptstyle \!\:\!t_{\mbox{}_{0}}}}^{\;\!t}
\!\!\:\!
\|\, D^{m + 1} \,\!\mbox{\boldmath $u$}(\cdot,\tau) \,
\|_{\mbox{}_{\scriptstyle L^{2}(\mathbb{R}^{4})}}
  ^{\:\!2}
d\tau
\;\leq\;\:\!
\|\, D^{m} \,\!\mbox{\boldmath $u$}(\cdot,t_0) \,
\|_{\mbox{}_{\scriptstyle L^{2}(\mathbb{R}^{4})}}
  ^{\:\!2}
} $ \\
\mbox{} \vspace{-0.550cm} \\
\mbox{} \hfill (2.23) \\
\mbox{} \vspace{-0.550cm} \\
\mbox{} \hspace{+0.500cm}
$ {\displaystyle
+\;\,
K_{\!\;\!m}
\:\!
\sum_{\ell\,=\,0}^{m}
\;\!
\int_{\mbox{}_{\scriptstyle \!\:\!t_{\mbox{}_{0}}}}^{\;\!t}
\!\!\;\!
\|\, D^{\:\!\ell} \,\!\mbox{\boldmath $u$}(\cdot,\tau) \,
\|_{\mbox{}_{\scriptstyle L^{4}(\mathbb{R}^{4})}}
\;\!
\|\, D^{\:\!m - \ell} \,\!\mbox{\boldmath $u$}(\cdot,\tau) \,
\|_{\mbox{}_{\scriptstyle L^{4}(\mathbb{R}^{4})}}
\;\!
\|\, D^{\:\!m + 1} \,\!\mbox{\boldmath $u$}(\cdot,\tau) \,
\|_{\mbox{}_{\scriptstyle L^{2}(\mathbb{R}^{4})}}
\;\!
d\tau
} $ \\
\mbox{} \vspace{+0.175cm} \\
for some constant $ K_{\!\;\!m} $
that depends on $m$ only.
Recalling that
$ {\displaystyle
\;\!
\|\, \mbox{u} \,
\|_{\mbox{}_{\scriptstyle L^{4}(\mathbb{R}^{4})}}
\!\;\!\leq\,
\|\, D \:\!\mbox{u} \,
\|_{\mbox{}_{\scriptstyle L^{2}(\mathbb{R}^{4})}}
} $ \linebreak
\mbox{} \vspace{-0.575cm} \\
(for arbitrary
$ \mbox{u} \in H^{1}(\mathbb{R}^{4}) $),
we get,
using (2.22) above, \\
\mbox{} \vspace{-0.550cm} \\
\begin{equation}
\tag{2.24}
\|\, D^{\:\!\ell} \,\!\mbox{u} \,
\|_{\mbox{}_{\scriptstyle L^{4}(\mathbb{R}^{4})}}
\;\!
\|\, D^{\:\!m - \ell} \,\!\mbox{u} \,
\|_{\mbox{}_{\scriptstyle L^{4}(\mathbb{R}^{4})}}
\:\!\leq\;
\|\, D \:\!\mbox{u} \,
\|_{\mbox{}_{\scriptstyle L^{2}(\mathbb{R}^{4})}}
\;\!
\|\, D^{\:\!m + 1} \,\!\mbox{u} \,
\|_{\mbox{}_{\scriptstyle L^{2}(\mathbb{R}^{4})}}
\end{equation}
\mbox{} \vspace{-0.100cm} \\
for every
$ \;\! 0 \leq \ell \leq m $,
so that
(2.23) gives \\
\mbox{} \vspace{+0.000cm} \\
\mbox{} \hspace{+0.750cm}
$ {\displaystyle
\|\, D^{m} \,\!\mbox{\boldmath $u$}(\cdot,t) \,
\|_{\mbox{}_{\scriptstyle L^{2}(\mathbb{R}^{4})}}
  ^{\:\!2}
\!\:\!+\;
2 \, \nu \!\!\;\!
\int_{\mbox{}_{\scriptstyle \!\:\!t_{\mbox{}_{0}}}}^{\;\!t}
\!\!\:\!
\|\, D^{m + 1} \,\!\mbox{\boldmath $u$}(\cdot,\tau) \,
\|_{\mbox{}_{\scriptstyle L^{2}(\mathbb{R}^{4})}}
  ^{\:\!2}
d\tau
\;\leq\;\:\!
\|\, D^{m} \,\!\mbox{\boldmath $u$}(\cdot,t_0) \,
\|_{\mbox{}_{\scriptstyle L^{2}(\mathbb{R}^{4})}}
  ^{\:\!2}
} $ \\
\mbox{} \vspace{-0.525cm} \\
\mbox{} \hfill (2.25) \\
\mbox{} \vspace{-0.500cm} \\
\mbox{} \hspace{+1.750cm}
$ {\displaystyle
+\;\;\!\;\!
(\:\!m + 1)
\,
K_{\!\;\!m}
\!\;\!
\int_{\mbox{}_{\scriptstyle \!\:\!t_{\mbox{}_{0}}}}^{\;\!t}
\!\!\;\!
\|\, D \:\!\mbox{\boldmath $u$}(\cdot,\tau) \,
\|_{\mbox{}_{\scriptstyle L^{2}(\mathbb{R}^{4})}}
\;\!
\|\, D^{\:\!m + 1} \,\!\mbox{\boldmath $u$}(\cdot,\tau) \,
\|_{\mbox{}_{\scriptstyle L^{2}(\mathbb{R}^{4})}}
  ^{\:\!2}
\,
d\tau
} $. \\
\mbox{} \vspace{+0.175cm} \\
In particular,
choosing
$ \;\! t_0 >\:\! t_{\ast} $
so that
$ {\displaystyle
\;\!
(\:\!m + 1) \;\!\;\! K_{\!\;\!m} \,
\|\, D\:\!\mbox{\boldmath $u$}(\cdot,\tau) \,
\|_{\mbox{}_{\scriptstyle L^{2}(\mathbb{R}^{4})}}
\!\;\!<\, \nu
\;\!
} $
for all
$ \:\!\tau > t_ 0 $
\mbox{[}\,which is possible
because
$ {\displaystyle
\;\!
\lim_{\tau \,\rightarrow\, \infty}
\|\, D \:\! \mbox{\boldmath $u$}(\cdot,\tau) \,
\|_{\mbox{}_{\scriptstyle L^{2}(\mathbb{R}^{4})}}
\!\;\!=\;\! 0
} $
(cf.\;(2.18)\,\!)\;\!\mbox{]},
it follows that \\
\mbox{} \vspace{-0.500cm} \\
\begin{equation}
\tag{2.26}
\|\, D^{m} \,\!\mbox{\boldmath $u$}(\cdot,t) \,
\|_{\mbox{}_{\scriptstyle L^{2}(\mathbb{R}^{4})}}
  ^{\:\!2}
\!\:\!+\;\;\!
\nu \!\!\;\!
\int_{\mbox{}_{\mbox{\footnotesize \!\:\!$s$}}}
    ^{\mbox{\footnotesize \;\!$t$}}
\!\!\;\!
\|\, D^{m + 1} \,\!\mbox{\boldmath $u$}(\cdot,\tau) \,
\|_{\mbox{}_{\scriptstyle L^{2}(\mathbb{R}^{4})}}
  ^{\:\!2}
d\tau
\;\leq\;\:\!
\|\, D^{m} \,\!\mbox{\boldmath $u$}(\cdot,s) \,
\|_{\mbox{}_{\scriptstyle L^{2}(\mathbb{R}^{4})}}
  ^{\:\!2}
\end{equation}
\mbox{} \vspace{-0.125cm} \\
for all
$ \;\! t > s > t_{0} $,
showing that
$ {\displaystyle
\;\!
\|\, D^{m} \,\!\mbox{\boldmath $u$}(\cdot,t) \,
\|_{\mbox{}_{\scriptstyle L^{2}(\mathbb{R}^{4})}}
\!\;\!
} $
is monotonically decreasing in
$ \;\![\;\!t_0, \:\!\infty ) $.
This concludes the proof
when $ \:\! n = 4 $.
The cases $ \:\!n = 2, \:\!3 $
are handled in a similar way,
using (2.19) instead of (2.18)
and appropriate
replacements for (2.23) and (2.24).
\mbox{} \hfill $ \Box $
}
%
% ------------------------------------------- End of Proof
%                                            for Theorem 2.3

%
% *************************************************************
% *                                                           *
% *                     Bibliography                          *
% *                                                           *
% *************************************************************
%
% *************************************************************
% *                                                           *
% *                      References                           *
% *                                                           *
% *************************************************************
%

%
% -------------------------------------------------------
%

%
\nl
\mbox{} \vspace{-0.250cm} \\
\nl
{\small
\begin{minipage}[t]{10.00cm}
\mbox{\normalsize \textsc{Pablo Gustavo Albuquerque Braz e Silva}} \\
Departmento de Matem\'atica \\
Universidade Federal de Pernambuco \\
Recife, PE 50740, Brazil \\
E-mail: {\sf pablo@dmat.ufpe.br} \\
\mbox{} \vspace{-0.620cm} \\
\mbox{} \hspace{+1.150cm}
        {\sf braz.pablo@gmail.com} \\
\end{minipage}
\nl
\mbox{} \vspace{-0.450cm} \\
\nl
\begin{minipage}[t]{10.00cm}
\mbox{\normalsize \textsc{Jana\'\i na Pires Zingano}} \\
Departamento de Matem\'atica Pura e Aplicada \\
Universidade Federal do Rio Grande do Sul \\
Porto Alegre, RS 91509, Brazil \\
E-mail: {\sf jzingano@mat.ufrgs.br} \\
\mbox{} \hspace{+1.150cm}
        {\sf janaina.zingano@ufrgs.br} \\
\end{minipage}
\nl
\mbox{} \vspace{-0.450cm} \\
\nl
\begin{minipage}[t]{10.00cm}
\mbox{\normalsize \textsc{Paulo Ricardo de Avila Zingano}} \\
Departamento de Matem\'atica Pura e Aplicada \\
Universidade Federal do Rio Grande do Sul \\
Porto Alegre, RS 91509, Brazil \\
E-mail: {\sf paulo.zingano@ufrgs.br} \\
\mbox{} \hspace{+1.150cm}
        {\sf zingano@gmail.com} \\
\end{minipage}
}
%
% -------------------------------------------------------------
%

\end{document}